\documentclass[final, 1p]{elsarticle}  

\usepackage{amssymb,latexsym}
\usepackage{amsmath}

\biboptions{sort&compress}

\journal{: \: arXiv}

\begin{document}

\newtheorem{conjec}{Conjecture}
\newtheorem{lema}{Lemma}
\newtheorem{teo}{Theorem}
\newproof{prova}{Proof}
\newtheorem{corol}{Corollary}

\begin{frontmatter}

\title{Counterexamples to the conjectured transcendence of $\,\sum{1/(n+\alpha)^{k}}$, its closed-form summation and extensions to polygamma functions and zeta series}

\author{F.~M.~S.~Lima}

\address{Institute of Physics, University of Bras\'{i}lia, P.O.~Box 04455, 70919-970, Bras\'{i}lia, DF, Brazil}


\ead{fabio@fis.unb.br}


\begin{abstract}
In a recent work, Gun and co-workers have proposed that $\,\sum_{n=-\infty}^{\infty}{(n+\alpha)^{-k}}\,$ is a transcendental number for all integer $\,k$, $k > 1$, and $\,\alpha \in \mathbb{Q} \backslash \mathbb{Z}$.  Here in this work, this proposition is shown to be \emph{false} whenever $\,k\,$ is odd and $\,\alpha\,$ is a half-integer. It is also shown that these are the only counterexamples, which allows for a correct reformulation of the original proposition. This leads to a theorem yielding a closed-form expression for the summation of that series, which determines its arithmetic nature.  The result is then extended to a sum of polygamma functions and some related zeta series. In view of the recurrent appearance of these series and functions in different areas of mathematics and applications, the closed-form results put forward here could well be included in modern computer algebra systems (CAS). 
\medskip
\end{abstract}

\begin{keyword}
{\small Transcendental numbers \sep Cotangent derivatives \sep Polygamma function \sep Riemann Zeta function \sep Bernoulli polynomials \sep Symbolic computations}

\MSC[2010] 11J81 \sep 11M41 \sep 33B15 \sep 33F10 \sep 68W30
\end{keyword}

\end{frontmatter}

\section{Introduction}

In a recent work, Gun, Murty and Rath (GMR) have presented the following proposition for a certain Dirichlet series (see Theorem~4.1 in Ref.~\cite{GMR}).

\begin{conjec} [GMR]
\label{conj:GMR}
\;  Let $\,\alpha \in \mathbb{Q} \backslash \mathbb{Z}\,$ and $\,k>1\,$ be an integer. Then $\:\sum_{n=-\infty}^{+\infty}{1/(n+\alpha)^k}\,$ is a transcendental number.
\end{conjec}

Soon after its publication, I have shown in Ref.~\cite{Lima} that the main result of Ref.~\cite{GMR} (namely, its Theorem~3.1) is \emph{incorrect}. So, by suspecting that the above series could converge to an algebraic number for some $\alpha \in \mathbb{Q} \backslash \mathbb{Z}$, I have considered the possibility of finding a counterexample. Though the short proof provided there in Ref.~\cite{GMR} is, at first sight, pretty convincing, after a few tests I have found a simple one:  the series is \emph{null} for $\,\alpha = \frac12\,$ and $\,k=3$. Therefore, the original GMR proposition is \emph{incorrect}.

Here in this work, it is shown that, for all odd values of $\,k$, $k>1$, and half-integer values of $\,\alpha$, both the Theorem~4.1 of Ref.~\cite{GMR} and its proof are \emph{incorrect}.  By noting that these values of $\,\alpha\,$ and $\,k\,$ encompass all possible counterexamples, I have succeeded in reformulating the GMR proposition in the form of a new theorem which, on taking into account a recent result by Cvijovi\'{c}~\cite{Cvijo2009}, yields a closed-form expression for $\:\sum_{n=-\infty}^{+\infty}{{\,1/(n+\alpha)^k}}$, $\,\alpha \in \mathbb{Q} \, \backslash \, \mathbb{Z}$, in the form of an algebraic multiple of $\,\pi^{\,k}$.  As a consequence, the arithmetic nature of that series, as well as $\, \psi_{k}(1-\alpha) \,+ (-1)^{k+1} \, \psi_{k}(\alpha)\,$ and some related families of zeta series is determined.

\section{Counterexamples to the GMR proposition}

By testing the validity of Conjecture~\ref{conj:GMR} for $k=3$ and some rational numbers $\alpha \in (0,1)$, I have found the following counterexample.

\begin{lema}[First counterexample]
\label{lem:k3a12}
\quad  The series
\begin{equation*}
\sum_{n=\,-\infty}^{+\infty}{\frac{1}{\left(n+\frac12\right)^3}}
\end{equation*}
converges to zero (hence an algebraic number).
\end{lema}

\begin{prova}
\quad  By writing the above series as the sum of two series, one for the non-negative values of $\,n\,$ and the other for the negative ones, which is a valid procedure since that series converges absolutely, one has
\begin{eqnarray}
\sum_{n=-\infty}^{+\infty}{\frac{1}{\left(n+\frac12\right)^3}} &=& \sum_{n=-\infty}^{-1}{\frac{1}{\left(n+\frac12\right)^3}} + \sum_{n=0}^{\infty}{\frac{1}{\left(n+\frac12\right)^3}} \nonumber \\
&=& \sum_{n=-\infty}^{-1}{\frac{2^3}{\left(2\,n+1\right)^3}} + \sum_{n=0}^{\infty}{\frac{2^3}{\left(2\,n+1\right)^3}} \, .
\end{eqnarray}
Now, by substituting $n=-m$ in the series for $n<0\,$ and $\,n=j-1\,$ in the series for $n \ge 0$, one finds
\begin{equation*}
\sum_{n=-\infty}^{+\infty}{\frac{1}{\left(n+\frac12\right)^3}} =  -\,8 \sum_{m=1}^{\infty}{\frac{1}{(2\,m-1)^3}} + 8 \sum_{j=1}^{\infty}{\frac{1}{(2\,j-1)^3}} = 0 .
\end{equation*}
\begin{flushright} $\Box$ \end{flushright}
\end{prova}

The existence of a counterexample to the GMR conjecture implies that its original statement (see Theorem~4.1~(2) of Ref.~\cite{GMR}) is \emph{false}.  In fact, by scrutinizing the proof given there in Ref.~\cite{GMR}, I have found some defective points. For a better discussion, let us reproduce it below.

\begin{prova}[GMR \emph{incorrect proof}]
\quad  We know that
\begin{equation}
\sum_{n=-\infty}^{+\infty}{\frac{1}{(n+\alpha)^k}} = \frac{1}{\alpha} + \frac{(-1)^k}{(k-1)!} \, D^{k-1}\left(\pi\,\cot{\pi z}\right) \left.\vphantom{\hbox{\large (}} \! \right|_{z=\alpha} \, ,
\label{eq:GMRerro}
\end{equation}
where $\,D:=\frac{d}{dz}$. It is a consequence of a result of Okada~\cite{Okada} that $D^{k-1}\left(\pi\,\cot{\pi z}\right) \left.\vphantom{\hbox{\large (}} \! \right|_{z=\alpha}$ is non-zero. But then it is $\pi^k$ times a non-zero linear combination of algebraic numbers of the form $\csc{(\pi\,\alpha)}$, $\cot{(\pi\,\alpha)}$. Thus we have the result. 
\begin{flushright} $\Box$ \end{flushright}
\end{prova}

There in the Okada's cited work, Ref.~\cite{Okada}, one finds, in its only theorem, the following (correct) linear independence result.

\begin{lema}[Okada's theorem]
\label{lem:okada}
\; Let $k$ and $q$ be integers with $k>0$ and $q>2$. Let $T$ be a set of ${\,\varphi{(q)}/2}$ representatives \emph{mod} $q$ such that the union $\{T,-T\}$ is a complete set of residues prime to $q$. Then the real numbers $D^{k-1}(\cot{\pi z})\vert_{z={a/q}}$, $a \in T$, are linearly independent over $\mathbb{Q}$.\footnote{Here, $\varphi{(q)}$ is the Euler totient function.}
\end{lema}

\begin{prova}
See Ref.~\cite{Okada} for a proof of this lemma based upon the partial fraction decomposition of $D^{k-1}\left(\cot{\pi z}\right)$, valid for all $z \, {\not \in} \, \mathbb{Z}$, as well as a theorem by Baker-Birch-Wirsing on cyclotomic polynomials.
\begin{flushright} $\Box$ \end{flushright}
\end{prova}

It is important to emphasize that Lemma~\ref{lem:okada} says nothing about $\,D^{k}(\cot{\pi z})\vert_{z={a/q}}\,$ when $\,q=2$. In other words, the linear independence over $\mathbb{Q}$ is not guaranteed for the half-integer values of $\,z$, which makes these values a source of potential counterexamples for the GMR proposition. Moreover, the proof of Okada's theorem is based upon the following partial fraction decomposition of his function $\,F_k(z) = \frac{k}{(-2\pi i)^k} \: D^{k-1}\left(\pi\,\cot{\pi z}\right)$, valid for all complex $z \, {\not \in} \, \mathbb{Z}$ (see Eq.~(1) of Ref.~\cite{Okada}):
\begin{equation*}
-\,\frac{k!}{(2\pi i)^k} \sum_{n=-\infty}^{+\infty}{\frac{1}{(n+z)^k}} = \frac{k}{(-2\pi i)^k} \: D^{k-1}\left(\pi\,\cot{\pi z}\right) ,
\end{equation*}
which, for positive integer values of $k$, simplifies to
\begin{equation}
\sum_{n=-\infty}^{+\infty}{\frac{1}{(n+z)^k}} = \frac{(-1)^{k-1}}{(k-1)!} \: D^{k-1}\left(\pi\,\cot{\pi z}\right) .
\label{eq:SkAlfa}
\end{equation}
This equation, by itself, shows that the expression for $\sum_{n=-\infty}^{+\infty}{{\,1/(n+\alpha)^k}}$ taken into account by Gun and co-workers in their proof in Ref.~\cite{GMR}, reproduced in Eq.~\eqref{eq:GMRerro} above, is \emph{incorrect}.


In our search for more counterexamples it is suitable to define the following family of real functions whose domain is $z \in \mathbb{R} \backslash \mathbb{Z}\,$:\,\footnote{Note that a division by zero occurs whenever $\,m \in \mathbb{Z}$. Also, $\lim_{z \rightarrow m} |S_k(z)| = \infty$.}
\begin{equation}
S_k(z) := \sum_{n=-\infty}^{+\infty}{\frac{1}{(n+z)^k}} \, ,
\label{eq:defS}
\end{equation}
where $k$ is an integer, $k>1$. All functions $S_k(z)$ share the following mathematical properties.

\begin{lema}[Properties of $\,S_k(z)\,$]
\label{lem:rules}
\;  The functions $\,S_k(z)$ defined in Eq.~\eqref{eq:defS} present the following properties, valid for all $\,z\,$ in its domain:
\begin{itemize}
\item[(i)] All functions $\,S_k(z)\,$ are periodic in $\,z$, with an unitary period;
\item[(ii)] For even values of $\:k$, $S_k(z) > 0$;
\item[(iii)] All functions $\,S_k(z)\,$ are differentiable;
\item[(iv)] For odd values of $\,k$, $k>1$, $S_k(z)$ is a continuous, strictly decreasing function in each real interval $\,(m,m+1)$, $ \forall \, m \in \mathbb{Z}$.
\end{itemize}
\end{lema}

\begin{prova}
\quad  Property (i) follows from the fact that, $\forall \, k \in \mathbb{Z}$, $k>1$, and for all real $z \, {\not \in} \, \mathbb{Z}$, $S_k(z+1) = \sum_{n=-\infty}^{+\infty}{{\,1/\left[\,n+(z+1)\right]^k}} = \sum_{m=-\infty}^{+\infty}{{\,1/\left(m+z\right)^k}} = S_k(z)\,$.  Property (ii) is a consequence of the fact that, for any \emph{even} $k$, $k \ge 2$, every term ${\,1/(n+z)^k}\,$ of the series that defines $S_k(z)$ is positive.  Property (iii) follows from a term-by-term differentiation of the series for $S_k(z)$ with respect to $z$, $z \in (0,1)$, which does not affect the convergence in this \emph{open} interval (see, e.g., Theor.~10.3.11 of Ref.~\cite{Bloch}).\footnote{This also follows from the representation of $\,S_k(z)\,$ as a sum of two polygamma functions $\psi_k(z)$, as will be established in Eq.~\eqref{eq:fim}, since each $\psi_k(z)$ is differentiable for all $x \in (0,1)$.} The periodicity of $\,S_k(z)\,$ in $z$ is then taken into account to extend its differentiability to all real $z \, {\not \in} \, \mathbb{Z}\,$ (i.e., all points of its domain). Property (iv) follows from a less direct argument. Firstly, from Eqs.~\eqref{eq:SkAlfa} and~\eqref{eq:defS} we deduce that, for any positive integer $p$,
\begin{equation}
S_{2p+2}(z) = \frac{-1}{(2p+1)!} \: D^{2 p +1}\left(\pi\,\cot{\pi z}\right)
\label{eq:Spar}
\end{equation}
and
\begin{equation}
\frac{d }{d z} \, S_{2p+1}(z) = \frac{1}{(2 p)!} \: D^{2p+1}\left(\pi\,\cot{\pi z}\right) .
\label{eq:der1}
\end{equation}
On isolating the cotangent derivative in Eq.~\eqref{eq:der1} and substituting it in Eq.~\eqref{eq:Spar}, one finds
\begin{equation}
\frac{d }{d z} \, S_{2p+1}(z) = -\,(2p+1) \;\, S_{2p+2}(z) \, .
\label{eq:der1new}
\end{equation}
By property (ii), $S_{2p+2}(z) > 0$, thus ${\,d S_{2p+1}/d z} < 0$ for all real $z \, {\not \in} \, \mathbb{Z}$. Then, $S_{2p+1}(z)$ is a strictly decreasing function in each interval over which $S_{2p+1}(z)$ is a continuous function, namely the disjoint intervals $(m,m+1)$, $m \in \mathbb{Z}$.
\begin{flushright} $\Box$ \end{flushright}
\end{prova}
\medskip

An immediate consequence of the periodicity of the functions $\,S_k(z)\,$ is the repetition of the null result established for $\,\alpha = \frac12\,$ in Lemma~\ref{lem:k3a12}. 

\begin{lema}[Counterexamples for $\,k = 3$]
\label{lem:k3aqq}
\; For every integer $m$, the series
\begin{equation*}
\sum_{n=-\infty}^{+\infty}{\frac{1}{\left( n+m+\frac12 \right)^3}}
\end{equation*}
is \emph{null}, hence an algebraic number.
\end{lema}

\begin{prova}
\;  Since $S_3(z)$ is a periodic function (with unitary period), as proved in Lemma~\ref{lem:rules}, and $\,S_3\left(\frac12\right) = 0$, as shown in Lemma~\ref{lem:k3a12}, then $S_3\left(m+\frac12\right) = 0$.
\begin{flushright} $\Box$ \end{flushright}
\end{prova}

In fact, it can be shown that the counterexamples pointed out in Lemmas~\ref{lem:k3a12} and \ref{lem:k3aqq} are particular cases of a more general set.

\begin{lema}[Counterexamples for odd $\,k$]
\label{lem:kOdd}
For every \emph{odd} integer $\,k$, $\,{k>1}$, and every $\,m \in \mathbb{Z}$, the series
\begin{equation*}
\sum_{n=-\infty}^{+\infty}{\frac{1}{\left( n+m+\frac12 \right)^k}}
\end{equation*}
is \emph{null}, hence an algebraic number.
\end{lema}

\begin{prova}
\;  The proof for $S_{k}(\frac12)$ (i.e., for $m=0$), valid for any odd integer $k$, $k>1$, is analogue to that developed in Lemma~\ref{lem:k3a12} for $k=3$.  By writing the corresponding series as the sum of two series, one for $n<0$ and the other for $n \ge 0$, one has
\begin{eqnarray}
S_{k}\!\left(\frac12\right)  &=& \sum_{n=-\infty}^{+\infty}{\frac{1}{\left(n+\frac12\right)^k}} = \sum_{n=-\infty}^{-1}{\frac{1}{\left(n+\frac12\right)^k}} + \sum_{n=0}^{\infty}{\frac{1}{\left(n+\frac12\right)^k}} \nonumber \\
&=& \sum_{n=-\infty}^{-1}{\frac{2^k}{(2\,n+1)^k}} + \sum_{n=0}^{\infty}{\frac{2^k}{(2\,n+1)^k}} .
\end{eqnarray}
By substituting $n=-m$ in the series for $n<0$ and $n=j-1$ in the series for $n \ge 0$, one finds
\begin{equation}
S_{k}\left(\frac12\right) =  -\,2^k \, \sum_{m=1}^{\infty}{\frac{1}{(2\,m-1)^k}} + 2^k \, \sum_{j=1}^{\infty}{\frac{1}{(2\,j-1)^k}} = 0 \, .
\end{equation}
The extension of this null result to other half-integer values of $z$ follows from the periodicity of $\,S_k(z)\,$ in $\,z$, see property (i) of Lemma~\ref{lem:rules}.
\begin{flushright} $\Box$ \end{flushright}
\end{prova}

With these counterexamples and properties in hands, we can reformulate the GMR original proposition. This is done in the next section.

\section{Closed-form summations and arithmetic nature of $S_k(\alpha)$}

Let us present and prove a theorem which corrects the GMR original proposition, determining the summation of the series $\,S_k(\alpha)$ for every $k \in \mathbb{Z}$, $k>1$, and $\,\alpha \in \mathbb{Q} \backslash \mathbb{Z}$.  The theorem makes use of $B_n(x)$, the Bernoulli polynomial of degree $\,n$, implicitly defined by
\begin{equation}
\frac{t\,e^{t\,x}}{e^t -1} = \sum_{n=0}^\infty{B_n(x) \, \frac{\,t^n}{n!}} \, , \quad 0 < |t| < 2 \pi \, .
\end{equation}
It also makes use of the `floor' function $\lfloor x \rfloor$, $x \in \mathbb{R}$, which returns the greatest integer less than or equal to $x$.

\begin{teo}[Main result]
\label{teo:Main}
\quad Let $\,k\,$ be an integer, $k>1$.  Then, for every $\:\alpha \in \mathbb{Q} \, \backslash \, \mathbb{Z}$, the following closed-form summation holds:
\begin{equation*}
\sum_{n=-\infty}^{+\infty}{\frac{1}{\left(n+\alpha \right)^k}} = (-1)^{\left\lfloor \frac{k+2}{2} \right\rfloor} ~ \frac{(2 \pi q)^k}{q \,k!} \: \sum_{s=1}^q{B_k\!\left( \frac{s}{q} \right) \, f_k(2 \pi s \: \widetilde{\alpha}\,)} \, ,
\end{equation*}
where $\, \widetilde{\alpha} := \alpha -\lfloor \alpha \rfloor \,$ is the representative of $\alpha$ in $(0,1)$, $\,f_k(\theta) = \cos{\theta}\,$ if $\,k\,$ is even and ~$\sin{\theta}\,$ if $\,k\,$ is odd,  and $\,q>1\,$ is the denominator of $\:\widetilde{\alpha}$.
\end{teo}

\begin{prova}
\;  From Eqs.~(\ref{eq:SkAlfa})~and~(\ref{eq:defS}), we know that, for any integer $k$, $k>1$, and every real $z \, {\not \in} \, \mathbb{Z}$
\begin{equation}
S_{k}(z) = \sum_{n=-\infty}^{+\infty}{\frac{1}{(n+z)^k}} = \frac{(-1)^{k-1}}{(k-1)!} \; \pi \, D^{\,k-1}\left(\cot{\pi z}\right) . 
\label{eq:GMRbom}
\end{equation}
From Ref.~\cite{Cvijo2009} (see the only theorem at p.~218), we know that, if $n$, $p$ and $q$ are positive integers and $\,p<q$, then
\begin{equation}
D^{2n-1}\left(\cot{\pi z}\right) \left.\vphantom{\hbox{\large (}} \! \right|_{z=p/q} = (-1)^n \, \frac{(2 \pi q)^{2n-1}}{n} \: \sum_{s=1}^q{B_{2n}(s/q) \, \cos{(2 \pi s \, p/q)}}
\label{eq:cvi1}
\end{equation}
and
\begin{equation}
D^{2n}\left(\cot{\pi z}\right) \left.\vphantom{\hbox{\large (}} \! \right|_{z=p/q} = (-1)^{n-1} \, \frac{2\,(2 \pi q)^{2n}}{2n+1} \: \sum_{s=1}^q{B_{2n+1}(s/q) \, \sin{(2 \pi s \, p/q)}} \, .
\label{eq:cvi2}
\end{equation}
If $\,\alpha \in \mathbb{Q}\,$ belongs to $(0,1)$, we simply take $p/q = \alpha$. Otherwise, we can make use of the periodicity of $S_k(\alpha)$, see Lemma~\ref{lem:rules}, to find a rational $\widetilde{\alpha} \in (0,1)$ such that
\begin{equation}
\sum_{n=-\infty}^{+\infty}{\frac{1}{\left(n+\alpha \right)^k}} = \sum_{n=-\infty}^{+\infty}{\frac{1}{\left(n +\widetilde{\alpha} \, \right)^k}} \, .
\end{equation}
This is accomplished by taking ~$p/q = \widetilde{\alpha} = \alpha -\lfloor \alpha \rfloor$. The closed-form summation of $S_k(\widetilde{\alpha})$ follows by substituting $2n = k$ in Eq.~\eqref{eq:cvi1} and $2n+1 = k$ in Eq.~\eqref{eq:cvi2}, and then putting the results in the derivative in Eq.~\eqref{eq:GMRbom}.
\begin{flushright} $\Box$ \end{flushright}
\end{prova}

From Okada's theorem, our Lemma~\ref{lem:okada}, one knows that $\,D^{k-1}\left(\cot{\pi z}\right) \left.\vphantom{\hbox{\large (}} \! \right|_{z=\alpha}$ is \emph{non-zero}, the only possible exceptions being the half-integer values of $\alpha$, for which $\widetilde{\alpha} = \frac12$. Let $\,B_n = B_n(0) =B_n(1)\,$ be the $n$-th Bernoulli number. From Eqs.~\eqref{eq:cvi1} and~\eqref{eq:cvi2}, one readily finds $\,D^{2n-1}\left(\cot{\pi z}\right) \left.\vphantom{\hbox{\large (}} \! \right|_{z=1/2} = (-1)^n \, (2\pi)^{2n-1} \, (2^{2n}-1) \, B_{2n}/n\:$ and $\:D^{2n}\left(\cot{\pi z}\right) \left.\vphantom{\hbox{\large (}} \! \right|_{z=1/2} = 0$.\footnote{Some well-known properties of the Bernoulli numbers were taken into account here, e.g. $\,B_n = 0\,$ for every odd values of $n$, $n>1$, and $B_1 = -\frac12$.} Since, for all positive integer $n$, $B_{2n} \ne 0$, then these derivatives can be taken into account to determine the arithmetic nature of $\,S_k(\alpha)$, as follows.

\begin{corol}[Arithmetic nature of $S_k(\alpha)\,$]
\label{cor:algnat}
\;  For any integer $k$, $k>1$, and every $\,\alpha \in \mathbb{Q} \, \backslash \, \mathbb{Z}$, the series
$\: \sum_{n=-\infty}^{+\infty}{1/\left(n+\alpha \right)^k} \:$
is either null or an algebraic multiple of $\,\pi^{\,k}$. Moreover, it is null ~\textbf{if and only if} ~$k\,$ is odd and $\,\alpha\,$ is a half-integer.
\end{corol}

\begin{prova}
As explained just above this corollary, the series $S_k(\alpha)$ is null if and only if $\,k\,$ is odd and $\,\widetilde{\alpha}=\frac12$. This conclusion can also be reached by noting that when $\,D^{k-1}\left(\cot{\pi z}\right) \left.\vphantom{\hbox{\large (}} \! \right|_{z=\alpha}=0$, the resulting equation $S_k(\alpha) = 0$ has no real roots if $k$ is even, according to property (ii) of Lemma~\ref{lem:rules}. For odd values of $k$, on the other hand, all half-integer values of $\,\alpha$ are roots of $S_k(\alpha) = 0$, as established in Lemma~\ref{lem:kOdd}.  All that rests is to show that the half-integer values of $\alpha$ are the \emph{only} roots of $S_{2n+1}(\alpha) = 0$. In fact, in the interval $(0,1)$, $\alpha = \frac12\,$ is a root of $S_{2n+1}(\alpha) = 0$ for all positive integer $n$, as guaranteed by Lemma~\ref{lem:kOdd}.  From properties (iii) and (iv) of Lemma~\ref{lem:rules}, we know that $S_{2n+1}(z)$ is a \emph{strictly decreasing} differentiable function, for all real $z \in (0,1)$. Therefore, $S_{2n+1}(z)$ cannot be null for two distinct values of $z$ within the interval $(0,1)$ and then $\,z = \frac12\,$ is the \emph{only} root in this interval. This null result can be readily extended to all half-integer values of $\alpha$ due to the periodicity of $S_k(\alpha)$. When $S_k(\alpha) \ne 0$, it is an algebraic multiple of $\,\pi^{\,k}\,$ because Eqs.~\eqref{eq:cvi1} and~\eqref{eq:cvi2} imply that, $\forall \, \alpha \in \mathbb{Q} \backslash \mathbb{Z}$, $\alpha$ being not a half-integer, $\,D^{k-1}\left(\cot{\pi z}\right) \left.\vphantom{\hbox{\large (}} \! \right|_{z=\alpha}$ evaluates to $\,\pi^{k}\,$ times a non-null linear combination, with rational coefficients, of numbers of the form $\sin{\theta}$, $\cos{\theta}$, $\theta$ being a rational multiple of $\,\pi$, and these numbers are known to be both algebraic (see Sec.~III.4 of Ref.~\cite{Niven}).
\begin{flushright} $\Box$ \end{flushright}
\end{prova}

Since $\,\pi\,$ is a transcendental number, as first proved by Lindemann (1882), then $\pi^{\,k}$ is also transcendental, $\forall \: k \in \mathbb{Z}$, $k \ne 0$. This allows us to reformulate the GMR original proposition.

\begin{corol}[Reformulation of GMR proposition]
\label{cor:transcend}
\;  Let $\alpha \in \mathbb{Q} \backslash \mathbb{Z}$ and $k>1$ be an integer. Then
$\: \sum_{n=-\infty}^{+\infty}{1/(n +\alpha)^k} \:$
is either null or a transcendental number. It is null if and only if ~$k\,$ is odd and $\,\alpha\,$ is a half-integer.
\end{corol}

Another consequence of Theorem~\ref{teo:Main} emerges when we write the cotangent derivatives in Eq.~\eqref{eq:SkAlfa} in terms of the polygamma function $\,\psi_k(z) := \psi^{(k)}(z) = {\, d^k \psi(z)/dz^k}$, where $\psi_0(z) = \psi(z) := \frac{d}{d z} \ln{\Gamma{(z)}}$ is the so-called digamma function, $\Gamma{(z)}$ being the Euler gamma function. From the reflection formula for $\psi_k{(z)}$ (see, e.g. Ref.~\cite{Kolbig} and Sec.~5.15 of Ref.~\cite{Nist}), namely
\begin{equation}
\psi_k(1-z) +(-1)^{k+1} \, \psi_k(z) = (-1)^k \: \pi \: D^k(\cot{\pi z}) \, ,
\end{equation}
valid for all integer $k \ge 0$ and $z \in \mathbb{C} \backslash \mathbb{Z}$, and taking into account Eqs.~(\ref{eq:SkAlfa}) and~(\ref{eq:defS}), one has
\begin{equation}
S_k(z) = \frac{\psi_{k-1}(1-z) + (-1)^k \, \psi_{k-1}(z)}{(k-1)!} \, ,
\label{eq:fim}
\end{equation}
which holds for any integer $k>1$ and $z \not \in \mathbb{Z}$. From Corollary~\ref{cor:algnat}, we know that, for any integer $k>1$ and every $\alpha \in \mathbb{Q} \backslash \mathbb{Z}$, $S_k(\alpha)$ is either null or an algebraic multiple of $\,\pi^k$. This implies that $\,k! \: S_{k+1}(\alpha) = \psi_{k}(1-\alpha) \,- (-1)^{k} \, \psi_{k}(\alpha)\,$ is either null or an algebraic multiple of $\,\pi^{k+1}$. On taking Eqs.~\eqref{eq:cvi1} and~\eqref{eq:cvi2} into account, we can be indeed more specific.

\begin{corol}[Reflection formula for polygamma function]
\label{cor:psi}
\; With the same notation of Theorem~\ref{teo:Main}, for every integer $\,k$, $k \ge 1$, and $\,\alpha \in \mathbb{Q} \, \backslash \, \mathbb{Z}$,
\begin{equation*}
\psi_{k}(1-\alpha) \,+ (-1)^{k+1} \, \psi_{k}(\alpha) = -\,{(-1)}^{\left\lfloor\frac{k+1}{2}\right\rfloor} ~ \frac{(2 \pi q)^{k+1}}{q\,(k+1)} \, \sum_{s=1}^q{B_{k+1}\!\left( \frac{s}{q} \right) \, f_{k+1}(2 \pi s \, \widetilde{\alpha}\,)} \, ,
\end{equation*}
the sum being null \,\textbf{if and only if} $\:\alpha\,$ is a half-integer and $\,k\,$ is even.
\end{corol}

Finally, let us determine the arithmetic nature of a family of zeta series related to the reflection formula for $\,\psi_k(\alpha)$. In Eq.~(7) of Sec.~1.41 of Ref.~\cite{Grad}, one finds the Taylor series expansion
\begin{equation}
z \, \cot{z} = \sum_{n=0}^{\infty}{(-1)^n \, \frac{2^{2n} \, B_{2n}}{(2n)!} \: z^{2n}} \, ,
\end{equation}
which converges for all $\,z \in \mathbb{R}$, $0<|z|<\pi$.  By substituting $\,z\,$ by $\,\pi z\,$ and using the Euler's formula, namely
\begin{equation}
\zeta{(2n)} =(-1)^{n-1} \, \frac{2^{2 n-1} \,B_{2 n}\,\pi^{2 n}}{(2 n)!} \, ,
\end{equation}
where $\,\zeta{(s)} := \sum_{n=1}^\infty{1/n^s}$, $s > 1$, is the Riemann zeta function, it is easy to deduce that
\begin{equation}
\pi \, \cot{(\pi z)} = \frac{1}{z} - 2 \sum_{n=1}^{\infty}{\zeta{(2n)} \: z^{2n-1}} \, ,
\label{eq:zetaPar}
\end{equation}
the series being convergent for all real $z$ with $0<|z|<1$.\footnote{For $z=0$, the zeta series is null, but Eq.~\eqref{eq:zetaPar} is not valid due to divisions by zero.  However, it remains valid for $z \rightarrow 0^{+}$, as the reader can easily check.}  By calculating successive derivatives on both sides of Eq.~\eqref{eq:zetaPar}, one easily obtains the following formulae for the cotangent derivatives of order $\,m$, $\,m \ge 0\,$:
\begin{equation}
D^{2m} (\pi \, \cot{\pi z}) = \frac{(2m)!}{z^{2m+1}} \,-2 \sum_{n=m+1}^{\infty}{\zeta{(2n)} \cdot (2n-1) \cdots (2n-2m) \; z^{2n-2m-1}}
\label{eq:derPar}
\end{equation}
and
\begin{equation}
D^{2m+1} (\pi \, \cot{\pi z}) = -\,\frac{(2m+1)!}{z^{2m+2}} \: ~ -2 \sum_{n=m+1}^{\infty}\!\zeta{(2n)} \cdot (2n-1) \cdots (2n-2m-1) \; z^{2n-2m-2}
\label{eq:derImpar}
\end{equation}
From Eqs.~\eqref{eq:cvi1} and~\eqref{eq:cvi2}, by substituting $\,n -m = j\,$ we readily find closed-form summations for these zeta series.


\begin{corol}[Related zeta series]
\label{cor:zeta}
\; With the same notation of Theorem~\ref{teo:Main}, for every positive integer $\,m\,$ and rational $\:\alpha \in (-1,1) \, \backslash \, \{0\}$, the following summations hold:
\begin{eqnarray*}
\sum_{j=1}^{\infty}{\frac{(2j+2m-1)!}{(2j-1)!} ~\zeta{(2j+2m)} ~\alpha^{2j-1}} = \frac{(2m)!}{2\,\alpha^{2m+1}} \, +(-1)^{m} \: \pi \, \frac{(2 \pi q)^{2m}}{2m+1}\nonumber \\
 \times \, \sum_{s=1}^q{B_{2m+1}(s/q) \, \sin{(2 \pi s \, \widetilde{\alpha}\,)}}
\end{eqnarray*}
and
\begin{eqnarray*}
\sum_{j=1}^{\infty}{\frac{(2j+2m-1)!}{(2j-2)!} ~\zeta{(2j+2m)} ~\alpha^{2j-2}} = -\,\frac{(2m+1)!}{2\,\alpha^{2m+2}} \,+(-1)^m \, \pi \, \frac{(2 \pi q)^{2m+1}}{2m+2} \nonumber \\
\times \, \sum_{s=1}^q{B_{2m+2}(s/q) \, \cos{(2 \pi s \, \widetilde{\alpha}\,)}} .
\end{eqnarray*}
\end{corol}

Note that both these zeta series are the sum of a rational number and an algebraic multiple of an integer power of $\pi$, except for the former series with $\,\alpha = \pm \, \frac12$, when it reduces to $\: \pm \: 2^{\,2 m} \, (2 m)!$, hence a non-null integer.  These new closed-form expressions for zeta series could well be included in systematic collections of zeta series, such as that by Srivastava and Choi~\cite{LTSri}.

In view of the recurrent appearance of the series and functions studied here in different areas of science, from number theory to statistics and mathematical physics, the implementation of the closed-form expressions presented in this paper in modern computer algebra systems, such as \emph{Mathematica} and \emph{Maple}, is worth of consideration.



\end{document}